\documentclass[twocolumn,amssymb,aps,nofootinbib]{revtex4}
\usepackage{times,mathptmx,amsthm,amsmath}
\usepackage{pstricks,pst-node,graphics}
\makeatletter
\renewcommand{\thesection}{\arabic{section}}
\renewcommand{\p@subsection}{}
\renewcommand{\thesubsection}{\arabic{section}.\arabic{subsection}}
\renewcommand{\p@subsubsection}{}

\makeatother

\newcommand{\odisk}{\pscircle[fillstyle=solid,fillcolor=black]}
\newcommand{\G}{{\cal G}}
\newcommand{\HH}{{\cal H}}
\newcommand{\M}{{\bf M}}

\theoremstyle{definition}

\psset{linewidth=.5pt}

\begin{document}
\title{A reciprocity theorem for domino tilings}
\author{James Propp ({\tt propp@math.wisc.edu})}
\thanks{Supported by grants from the National Science
Foundation and the National Security Agency.} 
\affiliation{Department of Mathematics \\
University of Wisconsin - Madison}
\date{March 31, 2001}
\begin{abstract}
\vspace{.1in}
\noindent
Let $T(m,n)$ denote the number of ways to tile
an $m$-by-$n$ rectangle with dominos.
For any fixed $m$, the numbers $T(m,n)$
satisfy a linear recurrence relation,
and so may be extrapolated to negative values of $n$;
these extrapolated values satisfy the relation
$T(m,-2-n)=\epsilon_{m,n}T(m,n)$,
where $\epsilon_{m,n}=-1$ if $m \equiv 2$ (mod 4) and $n$ is odd
and where $\epsilon_{m,n}=+1$ otherwise.
This is equivalent to a fact demonstrated by Stanley 
using algebraic methods.
Here I give a proof that provides, among other things,
a uniform combinatorial interpretation of $T(m,n)$
that applies regardless of the sign of $n$.
\vspace{.1in}
\end{abstract}
\maketitle

\section{Introduction}
\label{sec-intro}

It has undoubtedly been observed many times that the Fibonacci sequence,
when run backwards as well as forwards from its initial conditions,
yields a doubly-infinite sequence 
that is symmetrical modulo some minus-signs:
$$...,-8,5,-3,2,-1,1,0,1,1,2,3,5,8,...$$
We will see that this
is just one instance of a more general phenomenon
involving domino tilings of rectangles.
For positive integers $m$ and $n$, let $T(m,n)$ denote
the number of ways to cover an $m$-by-$n$ rectangle
with 1-by-2 rectangles (``dominos'')
with pairwise disjoint interiors.
When $m = 2$, the values $T(m,n)$ form the Fibonacci sequence;
more generally, for each fixed $m$, 
the values of $T(m,n)$ form a sequence satisfying 
a higher-order linear recurrence relation
(see section \ref{sec-recur}
for an explanation of why such a recurrence relation must exist).

These recurrences allow one to extrapolate
$T(m,n)$ to negative values of $n$ in a natural way,
for each fixed positive integer $m$.
There are in fact infinitely many recurrences one could employ;
but, as section \ref{sec-recur} explains,
the extrapolated values of $T(m,n)$ do not depend on 
which of the recurrence relations is used.

One finds that as long as $m$ is not congruent to 2 (mod 4),
the sequence of values $T(m,n)$ (with $n$ varying
over the positive and negative integers)
exhibits palindromicity:
$$T(m,-2-n)=T(m,n) \ \mbox{for all $n$}.$$
When $m$ is congruent to 2 (mod 4),
a modified version of this holds:
$$T(m,-2-n)=(-1)^n T(m,n) \ \mbox{for all $n$}.$$

These properties of the sequence of values $T(m,n)$
can also be expressed in terms of
the generating function $F_m(x) = \sum_{n=0}^{\infty} T(m,n) x^n$,
with $T(m,0)$ taken to be 1.
This generating function is a rational function of $x$,
as a consequence of the recurrence relation governing the values $T(m,n)$.
Denote this rational function by $P_m(x)/Q_m(x)$.
Note that the polynomial $Q_m(x)$ encodes
the recurrence relation that governs 
the sequence of coefficients $T(m,n)$,
at least eventually;
if the degree of $P_m(x)$ is less than
the degree of $Q_m(x)$,
then the recurrence relation applies 
to all the coefficients of the generating function,
not just ones with $n$ sufficiently large.

Stanley showed \cite{St} (see parts (d) and (e) of his main Theorem)
that for $m>1$,
the degree of $Q_m$ is 2 more than the degree of $P_m$,
$Q_m(1/x) = Q_m(x)$ when $m \equiv 2$ (mod 4)
and $Q_m(1/x) = -Q_m(x)$ when $m \not\equiv 2$ (mod 4).
These reciprocity relations,
first observed empirically by Klarner and Pollack \cite{KP},
can be shown to be equivalent 
to the claims made earlier about the extrapolated values
of $T(m,n)$.

Stanley's proof makes use of algebraic methods;
in particular, it relies upon formulas 
proved by Kasteleyn \cite{Ka} and Temperley and Fisher \cite{TF}
that express $T(m,n)$ in terms of complex roots of unity
(or, if one prefers, trigonometric functions).

I will present a combinatorial way of thinking about the
two-sided infinite sequence $(T(m,n))_{n=-\infty}^{\infty}$. 
All features of the reciprocity relation ---
including the factor $\epsilon_{m,n}$ ---
will have direct combinatorial interpretations. 
Moreover, the method used here can also be used
to demonstrate similar reciprocity phenomena
for three- and higher-dimensional analogues of domino tiling.
For instance, if $T(k,m,n)$ denotes the number of ways
to tile a $k$-by-$m$-by-$n$ box with 1-by-1-by-2 blocks
(in any orientation),
then $T(k,m,-2-n)=\pm T(k,m,n)$,
where $T(k,m,-2-n)$ is obtained by taking
the recurrence relation for $T(k,m,n)$
(with $k,m$ fixed and $n$ varying)
and running it backward.
Such higher-dimensional tiling problems
are not amenable to Stanley's method,
since the formulas of Kasteleyn and Temperley and Fisher
do not extend beyond the two-dimensional setting.

\section{Signed matchings}
\label{sec-match}

Instead of using domino-tilings, we will use an equivalent mathematical model:
perfect matchings of grid-graphs,
also known as dimer covers or dimer configurations\footnote{The difference
between the titles of \cite{St} and \cite{KP},
compared with the difference in the names of the journals
in which the two articles appear,
clarifies the difference between dimers and dominos: 
Studying dimers configurations is applied mathematics.
Studying tilings by dominos is not.}.
A perfect matching of a graph $\G$
(hereafter called just a matching for short)
is a subset of the edge-set of $\G$
such that each vertex of $\G$
belongs to one and only one edge of the subset.
The domino tilings of an $m$-by-$n$ rectangle
are in obvious bijection with
the perfect matchings of a graph with $mn$ vertices
arranged in $m$ rows of $n$ vertices each,
with edges joining horizontal and vertical neighbors.
See Figure 1.
$$
\begin{pspicture}(-2,0)(2,1)
\multips(-2,0)(1,0){5}{\odisk{.1}}
\multips(-2,1)(1,0){5}{\odisk{.1}}
\multips(-2,0)(1,0){4}{\psline(0,0)(1,0)}
\multips(-2,1)(1,0){4}{\psline(0,0)(1,0)}
\multips(-2,0)(1,0){5}{\psline(0,0)(0,1)}
\end{pspicture}
$$
\begin{center} 
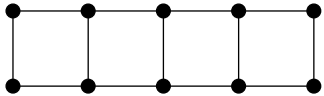
Fig. 1: The graph $\G(2,5)$.
\end{center}

Given that some of the extrapolated numbers $T(m,n)$ are negative,
it is natural to try to provide them with a combinatorial interpretation
by way of a signed version of graph-matching.
A signed graph is a graph whose edges are labeled $+1$ and $-1$,
and a matching counts as either positive or negative
according to whether the number of $-1$ edges is even or odd.
(This is a special case of weighted matching,
in which each edge is assigned a weight,
and the weight of a matching is the product of the weights of its edges.)
When we count the matchings of a signed graph,
we count them according to sign.
For example, the Fibonacci number $T(2,-5)=-3$
will arise as the ``number'' of signed matchings
of the signed graph shown in Figure 2.
$$
\begin{pspicture}(-3,0)(3,1)
\multips(-3,0)(1,0){7}{\odisk{.1}}
\multips(-3,1)(1,0){7}{\odisk{.1}}
\multips(-3,0)(1,0){6}{\psline(0,0)(1,0)}
\multips(-3,1)(1,0){6}{\psline(0,0)(1,0)}
\multips(-2,0)(1,0){5}{\psline[linestyle=dashed](0,0)(0,1)}
\end{pspicture}
$$
\begin{center} 
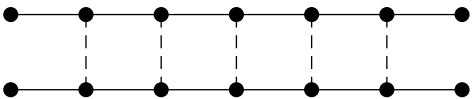
Fig. 2: The signed graph $\G(2,-5)$.
\end{center}
Here the horizontal edges are all $+1$ edges
and the vertical edges are $-1$ edges
(where the presence of a negative edge
is indicated by use of dashed lines).
Each of the three matchings of the graph
has an odd number of vertical edges
and so counts as a negative matching. 
Hereafter, when I refer to the number of matchings
of a signed graph,
I always mean the number of positive matchings
minus the number of negative matchings.

By a ``signed graph of width $m$'', I will mean
a subgraph of the $m$-by-$n$ grid graph for some $n$,
with some edges having sign 1 and others having sign $-1$.
We will be pay special attention to
some specific graphs $\G(m,n)$,
defined for $m \geq 1$ and $n$ an arbitrary integer.
When $n$ is positive, $\G(m,n)$ is just
the $m$-by-$n$ grid graph described above,
with all edges having sign 1.
When $n$ is less than or equal to zero,
$\G(m,n)$ is a modified version of the graph $\G(m,2-n)$
in which the rightmost and leftmost flanks of vertical edges
have been removed,
and all remaining vertical edges are decreed to have sign $-1$. 
Thus Figure 1 shows the graph $\G(2,5)$
while Figure 2 shows the graph $\G(2,-5)$.

It will emerge that the (signed) number of matchings of $\G(m,n)$
is $T(m,n)$.  
Note that for all $m$, $\G(m,0)$ and $\G(m,-2)$ each
have a unique matching, whose sign is $+1$. 
This verifies the reciprocity theorem
in the case $n=0$
(and the case $n=-1$ is trivial).
To prove the claim for other values of $n$
(it suffices to consider $n>0$),
we must do more work.

\section{Adjunction of graphs}
\label{sec-adjoin}

The two definitions of $\G(m,n)$ given above
(one for $n>0$, the other for $n \leq 0$)
are not just two distinct combinatorial notions, 
patched together at 0;
they fit ``seamlessly''.
To make this clearer, 
we will situate the graphs $\G(m,n)$ in a broader context,
in which they function as building blocks.
Suppose $\G_1$ and $\G_2$ are signed graphs of width $m$.
Define the ``adjoined graph'' $\G_1 \G_2$
as the signed graph we get when we place $\G_1$ to the left of $\G_2$
and join the rightmost vertices of $\G_1$
to the respective leftmost vertices of $\G_2$
using $n$ edges of type $+1$.
Figure 3 shows the graph $\G(2,3) \G(2,-3)$.
$$
\begin{pspicture}(1,0)(8,1)
\multips(1,0)(1,0){8}{\odisk{.1}}
\multips(1,1)(1,0){8}{\odisk{.1}}
\multips(1,0)(1,0){7}{\psline(0,0)(1,0)}
\multips(1,1)(1,0){7}{\psline(0,0)(1,0)}
\multips(1,0)(1,0){3}{\psline(0,0)(0,1)}
\multips(5,0)(1,0){3}{\psline[linestyle=dashed](0,0)(0,1)}
\end{pspicture}
$$
\begin{center} 
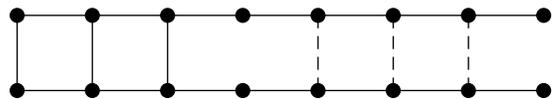
Fig. 3: The signed graph $\G(2,3) \G(2,-3)$.
\end{center}
Clearly adjunction is an associative operation,
so a product of three or more signed graphs of width $m$
is well-defined.

Let $M(\G)$ denote the (signed) number of matchings
of a signed graph $\G$.
I can now state the precise sense in which the two parts
of the definition of $\G(m,n)$ fit together:
$$
M(\G(m,n_1) \G(m,n_2) \cdots \G(m,n_k))
$$
\vspace{-0.3in}
$$
= M(\G(m,n_1+n_2+\dots+n_k))
$$
for all $n_1,\dots,n_k$, regardless of sign.

To indicate why the preceding equation is true,
we will prove the case $k=2$,
which contains all the essential ideas.

When $n_1$ and $n_2$ are positive,
we have $\G(m,n_1) \G(m,n_2)$ literally equal
to the graph $\G(m,n_1+n_2)$
by virtue of the definition of adjunction,
so this case requires no proof.

When $n_1$ and $n_2$ are both non-positive,
we get situations like the one shown in Figure 4,
depicting $\G(2,-2) \G(2,-3)$.
$$
\begin{pspicture}(-4,0)(4,1)
\multips(-4,0)(1,0){9}{\odisk{.1}}
\multips(-4,1)(1,0){9}{\odisk{.1}}
\multips(-4,0)(1,0){8}{\psline(0,0)(1,0)}
\multips(-4,1)(1,0){8}{\psline(0,0)(1,0)}
\multips(-3,0)(1,0){2}{\psline[linestyle=dashed](0,0)(0,1)}
\multips(1,0)(1,0){3}{\psline[linestyle=dashed](0,0)(0,1)}
\end{pspicture}
$$
\begin{center} 
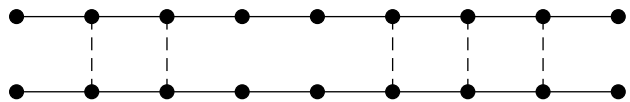
Fig. 4: The signed graph $\G(2,-2) \G(2,-3)$.
\end{center}
It is a well-known lemma in the theory of matchings
that the number of matchings of a graph $\G$
is unaffected if a chain of vertices $u,v,w$
(where $v$ is connected to $u$ and $w$ but
to no other vertices in the graph)
is shrunk down to a single vertex
(i.e., $v$ is removed and $u$ and $w$ are replaced
by a single vertex
adjacent to all the vertices other than $v$
to which $u$ and $w$ were originally connected).
The same observation applies to signed graphs as well
if the edges $uv$ and $vw$ are $+1$ edges.
Applying this lemma to $\G(2,-2) \G(2,-3)$ twice
(once in the top row and once in the bottom row),
we obtain a copy of the graph $\G(2,-5)$.
More generally, $M(\G(m,-a) \G(m,-b)) = M(\G(m,-(a+b)))$
whenever $a,b \geq 0$.

When $n_1$ and $n_2$ are of opposite type
(one positive, the other non-positive),
we get situations like the one shown in Figure 3,
depicting $\G(2,3) \G(2,-3)$.
If we apply the shrinking lemma in this case,
once in each row,
we get the signed graph shown in Figure 5.
$$
\begin{pspicture}(1,0)(6,1)
\multips(1,0)(1,0){6}{\odisk{.1}}
\multips(1,1)(1,0){6}{\odisk{.1}}
\multips(1,0)(1,0){5}{\psline(0,0)(1,0)}
\multips(1,1)(1,0){5}{\psline(0,0)(1,0)}
\multips(1,0)(1,0){2}{\psline(0,0)(0,1)}
\psline(3,0)(3,1)
\rput(2.8,0.5){$+$}
\rput(3.2,0.5){$-$}
\multips(4,0)(1,0){2}{\psline[linestyle=dashed](0,0)(0,1)}
\end{pspicture}
$$
\begin{center} 
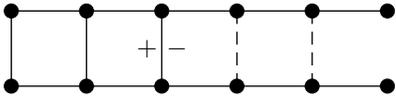
Fig. 5: $\G(2,3) \G(2,-3)$, after shrinking.
\end{center}
Here the edge marked $+-$ is actually two edges,
one a $+1$ edge and the other a $-1$ edge.
Any signed matching that uses the $+1$ edge
can be paired with a signed matching that uses the $-1$ edge;
these matchings have opposite sign,
and so together contribute 0 to the signed number of matchings.
Hence the number of matchings of the graph is unaffected
if we remove both of these edges,
obtaining the signed graph shown in Figure 6.
$$
\begin{pspicture}(1,0)(6,1)
\multips(1,0)(1,0){6}{\odisk{.1}}
\multips(1,1)(1,0){6}{\odisk{.1}}
\multips(1,0)(1,0){5}{\psline(0,0)(1,0)}
\multips(1,1)(1,0){5}{\psline(0,0)(1,0)}
\multips(1,0)(1,0){2}{\psline(0,0)(0,1)}
\multips(4,0)(1,0){2}{\psline[linestyle=dashed](0,0)(0,1)}
\end{pspicture}
$$
\begin{center} 
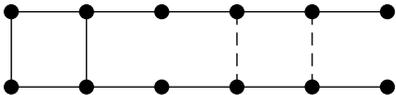
Fig. 6: $\G(2,3) \G(2,-3)$, reduced to $\G(2,2) \G(2,-2)$.
\end{center}
But this is just $\G(2,2) \G(2,-2)$.
More generally, applying the ``shrinking lemma'' to
$\G(m,a) \G(m,-b)$ yields $\G(m,a-1) \G(m,-(b-1))$
as long as $a,b$ are both positive.
This sets the stage for a proof by induction,
and all that remains is to verify the base case
where $a=0$ or $b=0$.
This is an easy verification that I leave to the reader.
The case of $\G(m,-a) \G(m,b)$ is identical.

\section{Recurrence relations}
\label{sec-recur}

We start by showing that
for every width-$m$ signed graph $\HH$,
the sequence $M(\HH \G(m,n))$ (with $n=1,2,3,\dots$)
satisfies a linear recurrence relation of degree at most $2^m$
that depends only on $m$, not on $\HH$.

Let $\HH$ be a width-$m$ signed graph.
For any sequence $c_1,\dots,c_m$ of $m$ bits,
and for $n \geq 1$,
let $S(n;c)$ be the sum of the signs of the matchings
of the signed graph
obtained from $\HH \G(m,n)$
by deleting the subset of the rightmost $m$ vertices
specified by the bit-pattern $c$;
specifically, we remove the rightmost vertex in row $i$
(and all edges incident with it)
for precisely those $i$ with $c_i=0$.
Let $S(n)$ be the row-vector of length $2^m$
whose components are the numbers $S(n;c)$.
Then (cf.\ \cite{KP} and \cite{Re})
there is a $2^m$-by-$2^m$ matrix $M$ such
that $S(n+1) = S(n) M$.
More specifically, the $c,c'$ entry of $M$
gives the number of ways in which 
each matching counted by $S(n;c)$
can be extended to a matching counted by $S(n+1;c')$,
adding only horizontal edges
between the $n$th and $n+1$st columns
and vertical edges in the $n+1$st column.
Hence we have $S(n) = S(1) M^{n-1}$.
By the Cayley-Hamilton Theorem,
the sequence of powers of $M$ satisfies the linear recurrence
associated with the characteristic polynomial of $M$,
and it follows that the sequence of vectors
$S(1), S(1) M, S(1) M^2, \dots$ also satisfies this recurrence.
In particular, the entry associated with $c=(1,1,\dots,1)$
satisfies this recurrence.

For instance, with $m=2$,
$M$ is the 4-by-4 matrix
$$
\left( \begin{array} {cccc}
0 & 0 & 0 & 1 \\
0 & 0 & 1 & 0 \\
0 & 1 & 0 & 0 \\
1 & 0 & 0 & 1 \end{array} \right),
$$
where rows and columns are indexed in the order
$(0,0),(0,1),(1,0),(1,1)$.
The characteristic polynomial of this matrix is
$x^4-x^3-2x^2+x+1$,
so the one-sided sequence of numbers
$a_n = M(\HH \G(m,n))$ (with $n$ going from 1 to $\infty$)
must satisfy the recurrence
$a_{n+4} - a_{n+3} - 2a_{n+2} + a_{n+1} + a_{n} = 0$,
regardless of the nature of the width-$m$ signed graph $\HH$.

I now claim that the two-sided sequence of numbers
$a_n$ (with $n$ going from $-\infty$ to $\infty$)
satisfies this same recurrence relation.
If we take $\HH=\G(m,-N)$ with $N$ large,
and use the fact that $M(\G(m,-N) \G(m,n)) = M(\G(m,n-N)) = T(m,n-N)$,
we see that if we take the doubly-infinite sequence 
$\dots,T(m,-1),T(m,0),T(m,1),\dots$
and start it at $T(m,-N)$,
we get a singly-infinite sequence
that (as we have shown) must satisfy
the Cayley-Hamilton recurrence.
Since $N$ was arbitrary,
the two-sided infinite sequence of $a_n$'s 
satisfies the Cayley-Hamilton recurrence as well.

Having constructed one recurrence relation
satisfied by the doubly-infinite sequence
${\bf v}=(\dots,T(m,-1),T(m,0),T(m,1),\dots)$,
or equivalently,
one linear operator $A$ that annihilates ${\bf v}$,
we must now consider others.
Could there be a linear operator $B$ that annihilates
the singly-infinite sequence 
$T(m,1),T(m,2),T(m,3),\dots$
but not the doubly-infinite sequence ${\bf v}$?
If so, then $B$ would yield a different extrapolation
of $T(m,n)$ to negative values of $n$
than $A$ would.

To show that no such operator $B$ exists,
note that $BA$ annihilates ${\bf v}$.
Since $BA=AB$,
$B$ must send ${\bf v}$ into
something annihilated by $A$. 
If $B$ does not itself annihilate ${\bf v}$,
then $B$ sends ${\bf v}$
into some other sequence
that vanishes for all sufficiently large indices
but does not vanish for all indices.
It is easy to show that such a sequence
cannot be annihilated by $A$
or indeed by any linear operator
with constant coefficients.

The above arguments show that
the numbers $T(m,n)$,
defined for non-positive values of $n$ as above,
are the unique way to extrapolate
so that a linear recurrence is satisfied.

\section{Reciprocity}
\label{sec-recip}
With all the ingredients in place,
we can now give a very simple explanation of
the reciprocity relation for domino tilings
of rectangles of fixed width.
For concreteness, we start with the special case $T(2,-5)$.
Referring to Figure 2,
we see that every matching of $\G(2,-5)$
must contain the two leftmost horizontal edges
and the two rightmost horizontal edges.
If we remove these edges from the graph
(and the vertices incident with those edges, 
and the edges incident with those vertices),
we get a copy of $\G(2,3)$
in which the vertical edges have sign 1.
Since every matching of $\G(2,3)$
has an odd number of vertical edges,
we see that $T(2,-5)=-T(2,3)$.

More generally, when $n$ is positive,
every matching of $\G(m,-2-n)$
must contain the $m$ leftmost horizontal edges
and the $m$ rightmost horizontal edges.
Removing these edges (and concomitant vertices and edges)
leaves a copy of $\G(m,n)$
in which the vertical edges have sign 1.
It is well-known (see e.g.\ \cite{Th})
that every domino-tiling of a rectangle
(or indeed more general regions)
can be obtained from every other
by means of moves in which
a 2-by-2 block of horizontal dominos
is rotated to give
a 2-by-2 block of vertical dominos,
or vice versa.
Since moves of this kind do not affect the parity
of the number of vertical dominos (or, in our terms,
the parity of the number of vertical edges in a matching),
we know that that all the signed matchings
of our modified version of $\G(m,n)$
carry the same sign;
it remains to determine what this sign is.

If one is willing to appeal to the invariance-of-parity result
mentioned in the preceding paragraph,
it is simple to evaluate the common sign of all the matchings
by considering the all-horizontal matching
or the all-vertical matching.
But in the interests of making the article self-contained,
I give a direct argument (of a fairly standard kind).
Say we have a matching of $\G(m,n)$
that involves $k_i$ vertical edges
joining row $i$ and row $i+1$,
for $i$ ranging from 1 to $m-1$.
When these edges and their vertices are removed from the graph,
along with all other edges joining row $i$ and row $i+1$,
our graph splits into two subgraphs,
each of which must have an even number of vertices
(since we are assuming that we have a matching). 
Hence $ni-k_i$ and $n(m-i)-k_i$ must be even.
Using the congruence $k_i \equiv ni$ (mod 2),
we find that $k_1+\dots+k_{m-1}$,
the number of vertical edges,
must be congruent to $n(1+2+\dots+(m-1))=n(m-1)m/2$ (mod 2).
If $m$ is congruent to 0 or 1 (mod 4),
$(m-1)m/2$ is even,
so the number of vertical edges is even.
If $m$ is congruent to 3 (mod 4),
then (in order for a matching to exist)
$n$ must be even,
so that $n(m-1)m/2$ is even.
However, if $m$ is congruent to 2 (mod 4),
then $(m-1)m/2$ is odd,
so that $n(m-1)m/2$ is even or odd
according to the parity of $n$.

Hence,
if we define $\epsilon_{m,n}$
to be $-1$ when $m$ is congruent to 2 (mod 4)
and $n$ is odd,
and $+1$ otherwise,
we have
$$T(m,-2-n)=\epsilon_{m,n}T(m,n)$$
as claimed.

\section{Motivation}
\label{sec-motive}

The definition of $\G(m,n)$ for $n<0$
arose not by ad hoc insight but by
a particular vision of combinatorics 
as a specialized form of algebra.
Under this vision,
a collection of combinatorial objects
ought to be represented by a multivariate polynomial
whose coefficients are all equal to 1,
where the individual terms represent the combinatorial objects themselves
in some fashion.

In the case of matchings of $\G(2,n)$ with $n>0$,
we assign a formal variable (or ``weight'') to each edge,
define the weight of a matching as the product
of the weights of its constituent edges,
and define the polynomial $\M(\G(m,n))$
as the sum of the weights of all the matchings.
This is a polynomial in which every coefficient equals 1,
and in which the constituent monomials
encode the respective matchings of the graph.
If we think of $\G(2,n)$ as being embedded in
$\G(2,n+1)$ as the induced subgraph on the leftmost $n$ columns,
we find that the edge-variables occurring in $\M(\G(2,n))$
form a subset of the edge-variables occurring in $\M(\G(2,n+1))$,
and that indeed, there is a recurrence expressing
$\M(\G(2,n+1))$ in terms of $\M(\G(2,n))$ and $\M(\G(2,n-1))$:
$$\M(\G(2,n+1)) = y_{n+1} \M(\G(2,n)) + w_n x_n \M(\G(2,n-1)),$$
where $y_{n+1}$ is the weight of the rightmost vertical edge of $\G(2,n+1)$
and $w_n$ and $x_n$ are the weights 
of the rightmost horizontal edges of $\G(2,n+1)$.

We can run this recurrence backward,
obtaining in succession some rational functions of the variables $w_n,x_n,y_n$
with $n$ negative.
These rational functions are in fact Laurent polynomials;
moreover, all the coefficients of these Laurent polynomials are $\pm 1$.
Just as the monomials of $\M(\G(2,n))$ with $n>0$
encode the matchings of $\G(2,n)$,
one might assume that the monomials in
the extrapolated Laurent polynomials ``$\M(\G(2,n))$'' with $n \leq 0$
encode combinatorial objects that will yield
a combinatorial interpretation of the 
extrapolated values $T(2,n)$ with $n \leq 0$.
The monomials lead one to matchings of the graphs $\G(2,n)$ with $n<0$,
and once one has this signed-matching model,
the idea is readily extended to $T(m,n)$ for all $m$.

It is worth remarking that the reciprocity theorem
satisfied by the Laurent polynomials $\M(\G(m,n))$
is a bit more complicated than the reciprocity theorem
satisfied by the numbers $T(m,n)$;
in addition to the sign-factor $\epsilon_{m,n}$,
there is also a monomial factor involving the variables.
A comprehensive reciprocity theorem
for domino tilings of rectangles
would take this into account.
I intend to write a longer article that treats these issues,
and also applies them to other models,
such as matchings of honeycomb hexagons (see \cite{CLP}).
I expect that algebraic recurrences 
for sequences of Laurent polynomials
will also shed light
on the issues raised in Problem 32 (the last problem)
in \cite{Pr}.

\end{document}